\documentclass{amsart}          

\usepackage{amsmath,amsthm}     
\usepackage{amssymb}            
\usepackage{euscript}           
\usepackage{enumerate,calc}
\usepackage[matrix,arrow,curve,frame]{xy}    

\xymatrixcolsep{1.9pc}                          
\xymatrixrowsep{1.9pc}
\newdir{ >}{{}*!/-5pt/\dir{>}}                  

\newtheorem{thm}[subsection]{Theorem}
\newtheorem{defn}[subsection]{Definition}
\newtheorem{prop}[subsection]{Proposition}

\newtheorem{cor}[subsection]{Corollary}
\newtheorem{lemma}[subsection]{Lemma}

\theoremstyle{definition}  
\newtheorem{example}[subsection]{Example}

\newtheorem{remark}[subsection]{Remark}

  {\end{list}}

\newcommand{\dfn}{\textbf} 

\newcommand{\mdfn}[1]{\dfn{\mathversion{bold}#1}} 

\newcommand{\Wedge}             {\vee}

\newcommand{\tens}              {\otimes}               
\newcommand{\iso}               {\cong}  

\newcommand{\cat}{\EuScript}    
\newcommand{\cA}{{\cat A}}      
\newcommand{\cC}{{\cat C}}
\newcommand{\cD}{{\cat D}}
\newcommand{\cE}{{\cat E}}

\newcommand{\cM}{{\cat M}}
\newcommand{\cN}{{\cat N}}
\newcommand{\cO}{{\cat O}}

\newcommand{\cQ}{{\cat Q}}

\newcommand{\cT}{{\cat T}}

\newcommand{\Top}{{\cat Top}}
\newcommand{\Gtop}{G\Top}
\newcommand{\Spectra}{{\cat Spectra}}
\newcommand{\Set}{{\cat Set}}
\newcommand{\cCat}{{\cat Cat}}
\newcommand{\sSet}{s{\cat Set}}

\newcommand{\Man}{{\cat Man}}

\newcommand{\sPre}{sPre}
\newcommand{\Pre}{Pre}

\newcommand{\field}[1]  {\mathbb #1} 
\newcommand{\A}         {\field A}
\newcommand{\R}         {\field R}

\newcommand{\C}         {\field C}

\DeclareMathOperator*{\colim}{colim}
\DeclareMathOperator*{\hocolim}{hocolim}

\DeclareMathOperator{\sk}{sk}
\DeclareMathOperator{\coeq}{coeq}

\DeclareMathOperator{\Fact}{Fact}
\DeclareMathOperator{\coFact}{coFact}
\DeclareMathOperator{\sRes}{sRes}
\DeclareMathOperator{\coRes}{coRes}
\DeclareMathOperator{\Gr}{Gr}

\newcommand{\ra}{\rightarrow}                   
\newcommand{\lra}{\longrightarrow}              
\newcommand{\lla}{\longleftarrow}               
\newcommand{\llra}[1]{\stackrel{#1}{\lra}}      
\newcommand{\llla}[1]{\stackrel{#1}{\lla}}      

\newcommand{\we}{\llra{\sim}}                   
\newcommand{\bwe}{\llla{\sim}}
\newcommand{\cof}{\rightarrowtail}              

\newcommand{\inc}{\hookrightarrow}              
\newcommand{\dbra}{\rightrightarrows}           
\newcommand{\ldbra}{\overrightarrow{\underrightarrow{\phantom\lra}}}

\newcommand{\Id}{Id}                            

\newcommand{\mM}{\underline{\cM}}

\newcommand{\Cech}{\v{C}ech\ }

\newcommand{\re}{Re}

\newcommand{\sing}{Sing}

\newcommand{\assign}{\mapsto}

\newcommand{\ovcat}{\downarrow}

\newcommand{\adjoint}{\rightleftarrows}

\newcommand{\bdd}[1]{\partial\Delta^{#1}}
\newcommand{\del}[1]{\Delta^{#1}}

\newcommand{\he}{\simeq}

\newcommand{\pt}{pt}

\newcommand{\Smk}{Sm_k}

\newcommand{\SmC}{Sm/\C}
\newcommand{\sshv}{s\text{Shv}}
\newcommand{\MV}{{\cat M}{\cat V}}
\newcommand{\ovcf}{(\cC\times\Delta \ovcat F)}

\newcommand{\fix}{\mbox{}\par\noindent}

\newcommand{\dcoprod}[1]{{\displaystyle\coprod_{#1}}}
\numberwithin{equation}{subsection}

\begin{document}

\title{Universal Homotopy Theories}

\author{Daniel Dugger}
\address{Department of Mathematics\\ Purdue University\\ West
Lafayette, IN 47907 } 

\email{ddugger@math.purdue.edu}

\begin{abstract}
Begin with a small category $\cC$.  The goal of this short note is to
point out that there is such a thing as a `universal model category
built from $\cC$'.  We describe applications of this to the study of
homotopy colimits, the Dwyer-Kan theory of framings, to sheaf theory,
and to the homotopy theory of schemes.
\end{abstract}

\maketitle

\tableofcontents

\section{Introduction}
Model categories were introduced by Quillen \cite{Q} to provide a
framework through which one could do homotopy theory in various
settings.  They have been astonishingly successful in this regard, and
in recent years one of the first things one does when studying any
homotopical situation is to try to set up a model structure.  The aim
of this paper is to introduce a new, but very basic tool into the
study of model categories.

Our main observation is that given any small category $\cC$ it is
possible to expand $\cC$ into a model category in a very generic way,
essentially by formally adding homotopy colimits.  In this way one
obtains a `universal model category built from $\cC$'.  There is an
accompanying procedure which imposes relations into a model category,
also in a certain universal sense.  These two fundamental techniques
are the subject of this paper.  Although they are very formal---as any
universal constructions would be---we hope to indicate that these
ideas can be useful, and have some relevance to quite disparate areas
of model category theory.

There are two general themes to single out regarding this
material:
\begin{enumerate}[(1)]
\item Universal model categories give a method for creating a homotopy
theory from scratch, based on a category of `generators' and a set of
`relations'.  On the one hand this is a procedure for building model
categories in order to study some known phenomenon: in fact our
original motivation was to `explain' (if that can be considered the
right word) Morel and Voevodsky's construction of a homotopy theory
for schemes \cite{MV}.  On the other hand, however, this can be a
technique for understanding a model category one already has, by
asking what kinds of objects and relations are needed to reconstruct
that homotopy theory from the ground up (see section \ref{se:locex}).
\item If one is trying to prove a theorem which should hold in all
model categories---a generic result about the behavior of certain
homotopy colimits, for example---then very often it suffices to prove
the result just in some `universal' case.  Universal model categories
enjoy several nice properties---they are simplicial, proper,
cofibrantly-generated, etc.---and so when working in the universal
case one has a wealth of tools at one's disposal which are not
available in general model categories.  This gives a technique for
proving theorems analagous to a standard trick in algebra, whereby one
proves a result for all rings by first reducing to a universal example
like a polynomial ring. 
\end{enumerate}

In this paper our goal is to document the basic results about
universal model categories, and to generally discuss theme \#1.  The
second theme makes a brief appearance in section (\ref{se:hocolim}),
but a thorough treatment will be postponed for a future paper.

\medskip

\noindent
With somewhat more detail, here is an outline of the paper:

\smallskip

In Section 2 we construct a model category $U\cC$ and explain in what
sense it is the universal model category built from $\cC$.  This
generalizes a construction from category theory in which one
`completes' a category by formally adjoining colimits: the model
category $U\cC$ is in some sense obtained by formally adding {\it
homotopy} colimits (\ref{se:cofrep}).  The analogs here are very
precise: the category-theoretic construction involves looking at
categories of diagrams with values in $\Set$, whereas our
homotopy-theoretic analog uses diagrams of {\it simplicial} sets.

\smallskip


Section 3 deals with connections between our universal model
categories and the Dwyer-Kan theory of cosimplicial resolutions.
These resolutions are a tool for studying `higher-homotopies' in model
categories, and are used for example to obtain explicit formulas for
homotopy colimits.  Our message in this section is this: a resolution is
simply a map from a universal model category.  We explain in
(\ref{se:hocolim}) how homotopy colimits can be studied by `lifting'
them to the universal examples $U\cC$.  These universal examples are
actually {\it simplicial} model categories, and so the theory of
homotopy colimits is in this way reduced to a case which is very well
understood.

\smallskip

The model categories $U\cC$ are a kind of free object, like the free
group generated by a set.  Just as in the algebraic setting, there
turns out to be a way of `imposing relations' in model categories:
This is the well-known process of {\it localization}, which is
reviewed in section 5.  It is natural to then ask what kinds of
homotopy theories can be described by generators and relations---that
is, by starting with a universal model category $U\cC$ and then
localizing.  This question is the subject of section 6.

There is a very broad and useful class of model categories called the
{\it combinatorial} ones, which have been introduced by Jeff Smith.
We announce in (\ref{th:present}) the result that every combinatorial
model category is equivalent to a localization of some $U\cC$; the
proof is too involved to include here, but is instead given in the
companion paper \cite{D3}.  One immediate consequence is that every
combinatorial model category is equivalent to a model category which
is simplicial and left proper (the `simplicial' part was proven under
slightly more restrictive hypotheses in \cite{D1}, using very
different methods.)
\smallskip

Finally, in sections 7 and 8 we deal with some elementary
applications.  The first of these is an interpretation of Jardine's
model category of simplicial presheaves \cite{J2}: we point out that
giving a Grothendieck topology on a category $\cC$ amounts to
specifying certain `homotopy-colimit' type relations, and studying the
sheaf theory of $\cC$ is precisely studying the model category one
obtains from $U\cC$ by imposing those relations.

The second application is to the Morel-Voevodsky homotopy theory of
schemes.  We show that their constructions are equivalent to starting
with some basic category of schemes $\cC$, forming the universal
homotopy theory built from $\cC$, and then imposing certain
geometrically-natural relations.  All of this is very formal, and our
point is precisely that it {\it is} formal.   
 

\bigskip

We close this introduction by giving a precise, but brief description
of these universal model categories (this is done with more detail in
section 2.)  To describe the results we need a preliminary definition.
Suppose that $\cM$ and $\cN$ are model categories equipped with
functors $r\colon \cC \ra \cM$ and $\gamma\colon \cC \ra \cN$, as
depicted below:
\[ \xymatrix{\cC \ar[r]^{r}\ar[rd]_{\gamma} & \cM \\
                & \cN.
}
\]
We define a \mdfn{factorization of $\gamma$ through $\cM$} to be the
following data:
\begin{enumerate}[(i)]
\item A Quillen pair $L:\cM \adjoint \cN: R$, together with
\item a natural weak equivalence $\eta\colon L \circ r \we \gamma$.
\end{enumerate}
The factorization will be denoted by the triple $(L,R,\eta)$.  In this
paper it will be useful to regard a Quillen pair $L:\cM \adjoint \cN:R$
as a map of model categories $\cM \ra \cN$, which perhaps makes the term
`factorization' seem more appropriate.

We also define the \mdfn{category of factorizations
$\Fact_\cM(\gamma)$}: its objects are triples $(L,R,\eta)$ as above,
and a map $(L,R,\eta) \ra (L',R',\eta')$ is a natural transformation
$L \ra L'$ making the following diagram commute:
\[ \xymatrix{LrX \ar[rr]\ar[dr]_{\eta} && L'rX\ar[dl]^{\eta'} \\
                       & \gamma X
}
\]
(Note that giving a natural transformation $L \ra L'$ is
equivalent---via adjointness---to giving a transformation $R'\ra R$,
or to giving two maps $L\ra L'$ and $R' \ra R$ which are compatible with
the adjunctions.  So we could have adopted a more symmetric definition
of $\Fact_\cM(\gamma)$, but it would be equivalent to the above).

Here is the basic result:

\begin{prop}
There exists a closed model category $U\cC$ together with a functor
$r\colon\cC \ra U\cC$, such that the following is true: any map
$\gamma\colon\cC \ra \cM$ from $\cC$ to a model category factors
through $U\cC$, and moreover the category of such factorizations is
contractible.
\end{prop}

\noindent
When dealing with universal constructions in ordinary category theory
one typically finds that the category of choices is a contractible
groupoid---this is what is usually meant by saying that something is
`unique up to unique isomorphism'.  When working in the homotopical
setting, where the maps of interest are weak equivalences rather than
isomorphisms, a category of choices will rarely be a groupoid.  The
key property of `homotopically universal' constructions is that the
category of choices is contractible.  We may interpret the above
proposition as saying that $U\cC$ is the {\it universal model
category} built from $\cC$.  (Of course referring to {\it the}
universal model category is somewhat inappropriate, but we will
continue this abuse of language throughout the paper).  An explicit
construction of $U\cC$ as a diagram category is given in section 2.

\subsection{Organization of the paper}
We have already given a rough outline, but there are a couple of other
points to make.  The reader should be warned that section 3, dealing
with framings, is somewhat technical and not strictly necessary for
the rest of the paper.  It follows section 2 because they are closely
related, but many readers will want to read section 2 and then skip
ahead to section 5 their first time through.

We also need to give a warning about the proofs, which in most cases
we have kept very short, only giving a general indication of what one
should do.  This is on the one hand because of the formal nature of
the results: once one decides on what the correct definitions and
theorems are, then the results almost prove themselves.  On the other
hand there is always a certain amount of unpleasant machinery to be
dealt with, but inflicting this on the reader would distract from the
essentially simple character of the results.  Proofs which are
decidedly non-trivial are generally postponed until the very last
section of the paper, which the reader can refer back to when
necessary.

\subsection{Notation and terminology}

Our conventions regarding model categories, framings, and other
elements of abstract homotopy theory generally follow those of
Hirschhorn \cite{H}. Hovey's book \cite{Ho} is also a good reference.
In particular, model categories are assumed to contain small limits
and colimits, and to have functorial factorizations.

Following \cite{Ho} we will define a \dfn{map} of model categories
$L\colon\cM \lra \cN$ to be a Quillen pair $L:\cM\adjoint \cN:R$.
That is, a Quillen pair will be regarded as a map of model categories
in the direction of the left adjoint.  The results in this paper will
make it clear why this seems justified.  We use $L^{cof}(X)$ to denote
an object obtained by taking a cofibrant replacement for $X$ and then
applying $L$ to it, and $R^{fib}(Y)$ denotes the result of replacing
$Y$ by a fibrant object and then applying $R$.

If $\cC$ is a category and $X$, $Y$ are objects, then $\cC(X,Y)$
denotes the set of maps from $X$ to $Y$.  If $\cM$ is a model
category then $\mM(X,Y)$ denotes a homotopy function complex from $X$
to $Y$.

Finally, we must say something about our conventions regarding
homotopy colimits.  To define these in general model categories, the
approach taken in both \cite{DHK} and \cite{H} is to chose a framing
on the model category and then to define homotopy colimits via an
explicit formula.  The subtlety is that this yields a construction
which is only homotopy invariant for diagrams of cofibrant objects.
For us, however, when we write `$\hocolim$' we will {\it always} mean
something which is homotopy invariant: so our hocolim functors are
defined by first taking a functorial cofibrant-replacement of every
object in our diagram, and only then using the explicit formulas given
in \cite{DHK} or \cite{H}.  

\section{Universal model categories}

In this section we introduce the construction of universal model
categories, and indicate their basic properties.  This generalizes a
standard construction in category theory, by which one formally adds
colimits to a category $\cC$ by passing to the category of diagrams
$\Set^{\cC^{op}}$.  Our universal model category $U\cC$ is simply the
diagram category $\sSet^{\cC^{op}}$ equipped with an appropriate model
structure---it may loosely be thought of as the result of formally
adding homotopy colimits to $\cC$ (see (\ref{se:cofrep})).
Proposition~\ref{pr:univ} explains the universal property this
construction satisfies.

\bigskip

It will be helpful if we review a basic result from category theory.
Recall that a category $\cC$ is called {\it co-complete} if every
diagram in $\cC$ indexed by a small category has a colimit.  Given a
small category $\cC$, there is a universal way of expanding it into
a co-complete category, by considering the category of presheaves.

Recall that a \dfn{presheaf} on $\cC$ is simply a functor
$F\colon\cC^{op} \ra \Set$, and a map of presheaves is just a natural
transformation.  We will use $\Pre(\cC)$ to denote the category of
presheaves on $\cC$---this is just another name for the category of
diagrams $\Set^{\cC^{op}}$.  There is a canonical functor $r\colon\cC
\ra \Pre(\cC)$ called the Yoneda embedding, which sends an object $X
\in \cC$ to the presheaf $\,rX\colon Z \assign \cC(Z,X)$.  The object
$rX$ is called the `presheaf represented by $X$'.

One has the following standard result (cf.  \cite[Proposition 1.45(i)]
{AR}):

\begin{prop}\label{pr:presh}\fix
\begin{enumerate}[(a)]
\item Any functor $\gamma\colon\cC \ra \cD$ from $\cC$ into a
co-complete category $\cD$ may be factored through a
colimit-preserving map $\re\colon \Pre(\cC) \ra \cD$:
\[
\xymatrix{ \cC \ar[r]^-r\ar[dr]_{\gamma} & \Pre(\cC) \ar@{.>}[d]^{\re} \\
                      & \cD.
}
\]
Moreover, the factorization is unique up to unique isomorphism.
\item 
The map $\re$ comes equipped with a right adjoint $\sing: \cD \ra
\Pre(\cC)$.  
\end{enumerate}
\end{prop}

The proof will be left to the reader, but the basic
fact which makes it work is the observation that every presheaf
$F$ may be canonically expressed as a colimit of representables.  One
looks at the overcategory $\cC\downarrow F$ determined by the Yoneda
embedding $\cC \inc \Pre(\cC)$, and there is a canonical diagram
$(\cC\downarrow F)\ra \Pre(\cC)$ which sends $[rX\ra F]$ to $rX$.  The
colimit of this diagram is precisely $F$, and we'll usually write this
as
\begin{equation}
\label{eq:colim}
F\iso \colim_{rX\ra F} rX.
\end{equation}
$F$ may be thought of as the `formal colimit' of this diagram of
representables.  The functor $\re$ in the above proposition is a
`realization' functor, which takes a formal colimit in $\Pre(\cC)$ and
then builds it in the category $\cD$.  It's adjoint is the `singular'
functor, defined so that $\sing\, X$ is the presheaf $c\mapsto
\cD(\gamma c,X)$.

\begin{example}\fix
\begin{enumerate}[(a)]
\item Consider the simplicial indexing category $\Delta$.  The
category $\Pre(\Delta)$ is just the category of simplicial sets, and
the above result tells us that simplicial sets are just `formal
colimits' built from the basic simplices. 

There is an obvious functor $\Delta \ra \Top$ which sends $[n]$ to the
topological simplex $\Delta^n$.  Since $\Top$ is co-complete, the
above result gives an adjoint pair $\re: \Pre(\Delta) \adjoint \Top
:\sing$.  Of course these are just the usual realization and singular
functors.

\item There is also an obvious map $\Delta \ra \cCat$ into the
category of small categories: it sends $[n]$ to the category $\{0\ra 1
\ra \cdots \ra n\}$.  Since $\cCat$ is co-complete we immediately get
functors $\re: \Pre(\Delta) \adjoint \cCat: \sing$.  The functor
$\sing$ may be identified with the usual nerve construction, and the
functor $\re$ is the usual way of obtaining a category from a
simplicial set.  (I learned this nice example from Tibor Beke).

\end{enumerate}
\end{example}

Now let $U\cC$ denote the category $\sPre(\cC)$ of simplicial
presheaves on $\cC$.  There is an obvious embedding $\Pre(\cC) \ra
\sPre(\cC)$ which sends any presheaf $F$ to the {\it discrete\/}
simplicial presheaf containing $F$ in every dimension (with
identity maps as faces and degeneracies).  Throughout this paper we
will implicitly identify $\Pre(\cC)$ with its image in $U\cC$.
Composing this embedding with the Yoneda map $\cC \ra \Pre(\cC)$ gives
an embedding $\cC \inc U\cC$ which we will also call $r$, by abuse of
notation.

$U\cC$ is just the diagram category $\sSet^{\cC^{op}}$, and so we
can give it a model structure by saying that a map $F \ra G$ is a
\begin{enumerate}[(a)]
\item weak equivalence if every $F(X) \ra G(X)$ is a weak equivalence
in $\sSet$;
\item fibration if every $F(X) \ra G(X)$ is a fibration in $\sSet$;
\item cofibration if it has the left-lifting-property with respect to
the trivial fibrations.
\end{enumerate}

\noindent
This is called the \dfn{Bousfield-Kan model structure} (see
\cite[p. 314]{BK}).  It is known to be cofibrantly-generated, proper,
simplicial, and to have a wealth of other nice properties: it inherits
essentially any nice property of $\sSet$.  The weak equivalences are
generally called \dfn{objectwise weak equivalences}, and likewise for
the fibrations.

\begin{prop} 
\label{pr:univ}
Any functor $\gamma\colon\cC \ra \cM$ from $\cC$ into a model category
$\cM$ may be `factored' through $U\cC$ in the sense that there is a
Quillen pair $\re: U\cC \adjoint \cM: \sing$ and a natural weak
equivalence $\eta\colon \re\circ r \we \gamma$:
\[ \xymatrix{ \cC \ar[r]^r \ar[dr]_\gamma 
             & U\cC \ar[d]^{\re} 
                   \ar@{}[dl] |<<<<{\Downarrow} \\
            {} & \cM.
}
\]
Moreover, the category of such factorizations (as defined in the
introduction) is contractible.
\end{prop}

\begin{proof}[Idea of proof]
Proposition~\ref{pr:presh} allows us to extend $\gamma$ to an adjoint
pair of categories $\Pre(\cC) \adjoint \cM$.  To extend this further
to $U\cC$ we must add a simplicial direction, and figure out what the
realization of objects like `$X\tens \del{n}$' should be for $X\in
\cC$.  This is accomplished by the theory of cosimplicial resolutions,
discussed in the next section.  The proof will be completed at that
time.
\end{proof}

Note that the representables $rX$ are always cofibrant in $U\cC$, and
therefore their images $\re(rX)$ are cofibrant in $\cM$.  This is why
we needed the natural transformation $\eta$, because the above
triangle won't be able to commute on the nose unless $\gamma$ actually
took its values in the cofibrant objects.

\begin{example}
\label{ex:upt}
If we take $\pt$ to be the trivial category with one object and an
identity map, then $U(\pt)$ is just the model category $\sSet$.  So
the homotopy theory of simplicial sets is just the universal homotopy
theory on a point.  This is really a silly statement, as simplicial
sets are in some sense built into the very fabric of what people have
decided they mean by a `homotopy theory'.  We will see a more interesting
statement along these lines in Example~\ref{ex:udelta}.
\end{example}

\begin{example}
Let $G$ be a finite group, and let $\Gtop$ denote the usual model
category of $G$-spaces.  Consider the orbit category $\cO_G$, which is
the full subcategory of $\Gtop$ whose objects are the orbits $G/H$.
The inclusion $\cO_G \inc \Gtop$ gives rise to a Quillen pair
\[ U(\cO_G) \adjoint \Gtop, \]
and it is easy to check that the singular functor associates to any
$G$-space $X$ the diagram of (singular complexes of) its fixed spaces
$G/H \mapsto X^{H}$.  It is a classical theorem of equivariant
topology that this Quillen pair is an equivalence---the homotopy
theory of $G$-spaces is the universal homotopy theory generated by the
orbit category $\cO_G$.
\end{example}

\subsection{Cofibrant replacement in \mdfn{$U\cC$}}
\label{se:cofrep}

We close this section with two generalizations of (\ref{eq:colim}),
which will explain in what sense every object of $U\cC$ is a formal
homotopy colimit of objects in $\cC$.  This is accomplished by writing
down two very convenient cofibrant-replacement functors in $U\cC$.
Knowing nice versions of cofibrant-replacement is often an important
point in dealings with these model categories.

\medskip

Let $F$ be an object in $\Pre(\cC)$.  Define
$\tilde{Q}F$ to be the simplicial presheaf whose $n$th level is
\[ (\tilde{Q}F)_n=\coprod_{rX_n \ra \cdots \ra rX_0 \ra F} (rX_n) \]
and whose face and degeneracy maps are the obvious candidates ($d_i$
means omit $X_i$, etc.)  In the language of \cite{BK}, this is the
\dfn{simplicial replacement} of the canonical diagram $(\cC\downarrow
F) \ra \Pre(\cC)$.  (Notice that $\tilde{Q}F$ is in some sense the
formal homotopy colimit of this diagram.)  Also note that there is a
natural map $\tilde{Q}F \ra \pi_0(\tilde{Q}F)$, and that the codomain
is just $F$ by (\ref{eq:colim}).

\begin{lemma}
\label{le:cofrep}
The simplicial presheaf $\tilde{Q}F$ is cofibrant, and the map
$\tilde{Q}F \ra F$ is a weak equivalence.
\end{lemma}

If $F_*$ is an arbitrary simplicial presheaf then applying the functor
$\tilde Q$ in every dimension gives a bisimplicial presheaf, and we
let $QF$ denote the diagonal.  Once again there is a natural map
$QF\ra F$.

\begin{prop}[Resolution by representables]
\label{pr:cofrep}
For any simplicial presheaf $F$ one has that $QF$ is cofibrant, and
the map $QF \ra F$ is a weak equivalence.
\end{prop}

\begin{proof}[Proofs]
Both the lemma and proposition are proven in section \ref{ad:cofrep}.
\end{proof}

Notice that $QF$ is a simplicial presheaf which in every dimension is
a coproduct of representables.  We can think of it as the
realization---or homotopy colimit---of the diagram of representables
\[ 
\cdots
 \coprod_{rX_1 \ra rX_0 \ra F_1} (rX_1) 
         \ldbra
 \dcoprod{rX_0 \ra F_0} (rX_0),
\]
or as the `formal' homotopy colimit of the same diagram back in $\cC$.
The above proposition tells us that every simplicial presheaf is
canonically a homotopy colimit of representables, which of course is
the direct analog of (\ref{eq:colim}).

\medskip

There is another cofibrant-replacement functor for $U\cC$ which looks
a little different from the one above, but is useful in some
settings.  Consider the functor $\cC \times \Delta \ra U\cC$ defined
by $A \times [n] \mapsto rA \tens \del{n}$.  For a simplicial
presheaf $F\in U\cC$ we may form the overcategory $(\cC \times \Delta
\ovcat F)$, whose objects correspond to the data $[A\times [n], rA\tens
\del{n} \ra F]$.  This category comes equipped with a canonical
functor $(\cC\times\Delta \ovcat F) \ra U\cC$ sending $[A\times [n],
rA\tens \del{n} \ra F]$ to $rA\tens \del{n}$, and the colimit of this
functor is easily seen to be $F$ itself.  The {\it homotopy\/} colimit
is called the \mdfn{canonical homotopy colimit of $F$ with respect to
$\cC$}, and will be denoted by \mdfn{$\hocolim \ovcf$}.  (To form this
homotopy colimit recall that $U\cC$ is a simplicial model
category, and so we can use the formulas from \cite{BK}.)

The object $\hocolim \ovcf$ is a homotopy colimit of a diagram in
which all the objects have the form $rA\tens \del{n}$, and in
particular are cofibrant.  Therefore $\hocolim \ovcf$ is cofibrant as
well.  The natural map from the homotopy colimit to the colimit gives
a map $\hocolim \ovcf \ra F$, and we claim that this is always a weak
equivalence:

\begin{prop}
\label{pr:cofrep2}
Let $\cQ$ be the functor defined by $\cQ F=\hocolim \ovcf$.  Then
each $\cQ F$ is cofibrant, and the natural map $\cQ F \ra F$ is a weak
equivalence.  
\end{prop}

\begin{proof}
See section 9.
\end{proof}

The above proposition is of course another generalization of
(\ref{eq:colim}), as it shows how to canonically express any $F\in
U\cC$ as a homotopy colimit of representables.

\begin{example}
To see the difference between the functors $\cQ$ and $Q$ consider the
case where $\cC$ is the trivial category with one object and an
identity map.  Then $U\cC$ is the category $\sSet$.  Given a
simplicial set $K$, $Q K$ is just $K$ again.  To get $\cQ K$,
though, we take the category of simplices $\Delta K$ of $K$ and
consider the diagram $\Delta K \ra \sSet$ which sends the $n$-simplex
$\sigma$ to $\del{n}$.  The homotopy colimit of this diagram is
$\cQ K$.  It is weakly equivalent to $K$, but is a much bigger
object.
\end{example}

Canonical homotopy colimits are extremely important in \cite{D3},
where a treatment is given for all model categories.

\section{Connections with the theory of framings}

This section continues our discussion of the basic theory of universal
model categories.  What we will see is that working with these
universal model categories is exactly the same as working with {\it
cosimplicial resolutions\/}, in the sense of Dwyer and Kan.  Our
`universal' perspective has some advantages, however, in that it
efficiently captures the limited amount of adjointness that
resolutions exhibit.  We explain in (\ref{se:hocolim}) a technique by
which many theorems whose proofs require resolutions can be immediately
reduced to the case of {\it simplicial\/} model categories, which are
usually easier to deal with.  For instance, most standard results in
the theory of homotopy colimits can be deduced from the simplicial
case by this method.

\bigskip

We begin by reviewing what resolutions are.  First recall the notion
of a cylinder object: If $X\in \cM$ then a cylinder object for $X$ is an
object of $\cM$ which `looks and feels' like `$X\times \Delta^1$'.  It
is an object $X_1$ together with maps
\[ X \amalg X \cof X_1 \we X, \]
where the first map is a cofibration and the second a weak
equivalence.
These maps can be assembled into the beginning of a cosimplicial
object:
\[ \xymatrix{ X \ar[r]<0.5ex>\ar[r]<-0.5ex> & X_1. \ar@/_2ex/[l]
}
\]
The Dwyer-Kan theory of resolutions \cite{DK} is a massive
generalization of this, which gives a way of talking about objects
which `look and feel' like `$X \times \Delta^n$' for any $n$.  This is
actually what is called a {\it cosimplicial\/} resolution.  There are
also {\it simplicial\/} resolutions, which give a way of dealing with
objects which look and feel like $X^{\Delta^n}$, in the same way that
path objects are substitutes for $X^{\Delta^1}$.  The theories of
cosimplicial and simplicial resolutions are completely dual.

In a simplicial model category $\cM$ the object $X\tens \Delta^*$ is a
particularly nice element of $c\cM$---the category of cosimplicial
objects---with the property that the object in each level is weakly
equivalent to $X$ (at least if $X$ is cofibrant!)  The main part of
what one must come to terms with is what should be meant by
`particularly nice'---just as for cylinder objects, this should
translate into certain maps being cofibrations.  The reader can consult
\cite[Section 4.3]{DK} for a precise formulation.
It turns out that there is a natural model structure on $c\cM$ called
the {\it Reedy model structure\/} whose cofibrant objects are
precisely what we want.  We will not recall Reedy model categories
here, but refer the reader to \cite[Chapter 5]{Ho}.

\begin{defn} Let $\cM$ be a model category and let $X$ be an object.
\begin{enumerate}[(a)]
\item 
A \dfn{cosimplicial resolution} of $X$ is a Reedy cofibrant object
$\Gamma\in c\cM$ together with a map $\Gamma \ra c^*X$ which is a weak
equivalence in every dimension.  Here $c^*X$ is the constant
cosimplicial object with $X$ in every level.
\item
A \dfn{simplicial resolution} of $X$ is a Reedy fibrant object
$\Phi$ in $s\cM$ together with a map $c_*X\ra \Phi$ which is a weak
equivalence in every dimension.
\end{enumerate}
\end{defn}

\noindent
If $X$ is itself cofibrant then it has a cosimplicial resolution
whose $0$th object is actually equal to $X$, and the standard
practice is to choose such a resolution when possible.  

The above definition has the following immediate generalization:

\begin{defn}
Let $\cC$ be a category with a functor $\gamma\colon \cC \ra \cM$.  
A \dfn{cosimplicial resolution} of $\gamma$ is 
\begin{enumerate}[(i)]
\item a functor $\Gamma:\cC \ra c\cM$ such that each $\Gamma(X)$ is Reedy 
cofibrant, and
\item a natural weak equivalence $\Gamma(X) \we c^*X$,
\end{enumerate}

\dfn{Simplicial resolutions} of $\gamma$ are defined similarly.
\end{defn}

The convention in \cite{DHK} and \cite{H} is to use `framings' rather
than resolutions---the difference is, for instance, that in a
cosimplicial framing the objects $\Gamma(X)$ are Reedy cofibrant only
if $X$ itself was cofibrant.  The advantage of framings is that in a
simplicial model category the assignment $X \mapsto X\tens \Delta^*$
is a cosimplicial framing, whereas to get a cosimplicial resolution we
must use $X \mapsto X^{cof}\tens \Delta^*$ for some cofibrant
replacement $X^{cof} \we X$.  The disadvantage of framings is that
formulas which use them don't always yield the `correct' answer---to
get the correct answer one must use resolutions.

We will also need to talk about maps between resolutions:

\begin{defn}
Let $\Gamma_1$ and $\Gamma_2$ be two cosimplicial resolutions of a map
$\gamma:\cC \ra \cM$.  A \dfn{map of resolutions} $\Gamma_1 \ra
\Gamma_2$ is a natural transformation $\Gamma_1(X) \ra \Gamma_2(X)$
making the following triangle commute:
\[ \xymatrix{ \Gamma_1(X) \ar[rr]\ar[dr] && \Gamma_2(X) \ar[dl] \\
                      &c^*X }
\]
The category of cosimplicial resolutions on $\gamma$ will be denoted
\mdfn{$\coRes(\gamma)$}.
A map of simplicial resolutions, as well as the category
$\sRes(\gamma)$, are defined similarly.
\end{defn}

\begin{prop}
\label{pr:fact}
Let $\cC$ be a small category and let $\gamma\colon\cC \ra \cM$ be a map.
Then giving a factorization of $\gamma$ through $U\cC$ is precisely the
same as giving a cosimplicial resolution on $\gamma$.  Even more, there
is a natural equivalence of categories
\[ \Fact_\cM(\gamma) \he \coRes(\gamma).\] 
\end{prop}

\noindent
In words, the proposition says that a cosimplicial resolution is just a
map from a universal model category.

\begin{proof}
This is not hard, but requires some machinery. See section
\ref{ad:fact}.
\end{proof}

\begin{proof}[Proof of Proposition~\ref{pr:univ}]
We have just seen that factoring a functor $\gamma\colon\cC \ra \cM$
through $U\cC$ is equivalent to giving a cosimplicial resolution on
$\gamma$.  But it is a standard result in the theory of resolutions
that (1) any diagram $\gamma\colon \cC \ra \cM$ has a cosimplicial
resolution, and (2) the category $\coRes(\gamma)$ is contractible
(both are proven in \cite{H}).  So Proposition~\ref{pr:univ} is just a
re-casting of these classical facts.
\end{proof}


\subsection{Application to homotopy colimits}
\label{se:hocolim}\fix

Let $X\colon \cC \ra \cM$ be a diagram whose homotopy colimit we wish to
study.  When $\cM$ is a {\it simplicial\/} model category then there is
a well known formula for the homotopy colimit due to Bousfield and
Kan.  We can use universal model categories to reduce to the
simplicial case in a very natural way:

The map $X\colon\cC \ra \cM$ will factor through the universal model category
$U\cC$,
\[ \xymatrix {\cC \ar[r]^{\tilde X} \ar[dr]_{X} 
                        & U\cC \ar[d]^{\re} \\
                 & \cM.
}
\]
Note that $\tilde X$ doesn't really have anything to do with $X$---it
is the same functor we have been calling $r$, and it's the universal
example of a diagram in a model category with indexing category $\cC$.
For present purposes it's convenient to think of it as a lifting of
the diagram $X$, though.

$U\cC$ happens to be a simplicial model category, and so we can
use the Bousfield-Kan formula to compute the homotopy colimit of
$\tilde X$:
\[ \hocolim \tilde X \he \coeq \Biggl [ 
                          \coprod_{\beta \ra \gamma} \tilde X_\beta
                            \tens B(\gamma\downarrow \cC)^{op} 
                            \dbra 
                           \coprod_\alpha \tilde X_\alpha \tens
                                  B(\alpha\downarrow \cC)^{op}
                            \Biggr ].
\] 
Now the realization $\re$ is a left Quillen functor, and so whatever
we mean by the homotopy colimit of $X$ in $\cM$ will have to be weakly
equivalent to $\re (\hocolim \tilde X$).  
So we have uncovered the formula
\begin{align*}
 \hocolim X &\he \re (\hocolim \tilde X) \\
            & \he 
                  \coeq \Biggl [ 
                          \coprod_{\beta \ra \gamma} \re \bigl (\tilde X_\beta
                            \tens B(\gamma\downarrow \cC)^{op} \bigr ) 
                            \dbra 
                           \coprod_\alpha \re\bigl (\tilde X_\alpha \tens
                                  B(\alpha\downarrow \cC)^{op}\bigr )
                            \Biggr ]
\end{align*}

The Dwyer-Kan theory gives a formula for the homotopy colimit in terms
of resolutions: one chooses a cosimplicial resolution for the map $\cC
\ra \cM$, and then writes down an analog of the Bousfield-Kan formula
in which `tensoring' has been replaced by an operation involving
resolutions.  Our choice for the factorization of $X$ through $U\cC$
was precisely a choice of resolution for $X$, and it is an easy
exercise to check that the formula we obtained is precisely the
Dwyer-Kan formula.

One consequence of this perspective is that the basic properties of
homotopy colimits, once known for simplicial model categories,
immediately follow in the general case: one knows the property for the
universal case $U\cC$, and then simply pushes it to general $\cM$
using the left Quillen functor $\re$.  (The full power of this
technique requires the ability to impose `relations' into $U\cC$, as
discussed in section 5).  The idea is analagous to a standard trick in
algebra, where one proves a result for all rings by first reducing to
some universal example like a polynomial ring.  We hope to give a more
detailed treatment of this material in a future paper.

\section{Universal model categories for homotopy limits}

The categories $U\cC$ we have been talking about perhaps should have
been called `co-universal' model categories.  There is of course a
strictly dual notion which we will denote $V\cC$---this will be
discussed briefly in this section.  Just as $U\cC$ was very relevant
to the study of homotopy colimits, $V\cC$ pertains to the theory of
homotopy {\it limits\/}.  The material in this section is not needed
in the rest of the paper, but is included for the sake of
completeness.

\bigskip

Let $V\cC$ denote the category $[\sSet^{\cC}]^{op}$.  Note that there is
an obvious `Yoneda embedding' $r\colon\cC \ra V\cC$.  The diagram
category $\sSet^\cC$ may be given the usual Bousfield-Kan structure,
and we give $[\sSet^{\cC}]^{op}$ the opposite model structure: a map
$D_1 \ra D_2$ in $V\cC$ is a weak equivalence (resp. cofibration,
fibration) precisely when $D_2 \ra D_1$ is a weak equivalence
(resp. fibration, cofibration) in the model category $\sSet^{\cC}$.

Given two functors $\gamma\colon\cC \ra \cM$ and $r\colon\cC \ra \cN$
from $\cC$ to model categories $\cM$ and $\cN$, define a
\dfn{co-factorization of $\gamma$ through $\cN$} to be:
\begin{enumerate}[(i)]
\item A Quillen pair $L:\cM \adjoint \cN: R$, together with
\item a natural weak equivalence $\gamma \we R\circ r$.
\end{enumerate}
We leave it to the reader to define the \mdfn{category of
co-factorizations $\coFact_\cN(\gamma)$}: it is strictly dual to the
category of factorizations.  (Note that in this section it is more
convenient to think of a Quillen pair as a map of model categories in
the direction of the {\it right\/} adjoint!)

\begin{prop}\fix
\begin{enumerate}[(a)]
\item If $\gamma\colon\cC \ra \cM$ is a map from $\cC$ to a model
category then there is a co-factorization of $\gamma$ through $V\cC$,
and the category of all such co-factorizations is contractible.
\item The category of co-factorizations is naturally equivalent to the
category of simplicial resolutions on $\cC \ra \cM$.  
\end{enumerate}
\end{prop}

All of the results from sections 2 and 3 can be repeated in this
context.  In particular, the objects of $V\cC$ can be thought of as
formal homotopy limits of the objects of $\cC$.

\section{Imposing relations via localizaton}

Now that we have a notion of universal object for model categories, it
is natural to ask if there is some procedure for `imposing relations',
and then if every model category can be obtained from a universal one
in this way.  These questions will be addressed in this section and
the next.  Our method of imposing relations is the well-known
procedure of {\it localization\/}: given a model category $\cM$ and a
set of maps $S$, one forms a new model structure $\cM/S$ in which the
elements of $S$ have been added to the weak equivalences.  A very
thorough account of localization machinery is contained in \cite{H},
but in the beginning of this section we summarize the relevant
material.

Model categories of the form $U\cC/S$ will be our central concern for
the rest of the paper, and in (\ref{se:locex}) we give some basic
examples: the most notable of these is Segal's $\Gamma$-spaces, which
we can interpret in terms of universal constructions.  In
(\ref{se:cofibmc}) we end with some indications that the objects
$U\cC/S$ are something like `cofibrant' model categories.

\bigskip

\subsection{Review of localization}

Our basic definition of localizations for model categories is a slight
variant of what is called `left localization' by Hirschhorn \cite{H}:

\begin{defn} Let $\cM$ be a model category and let $S$ be a set of
maps in $\cM$.  An \mdfn{$S$-localization of $\cM$} is a model
category $\cM/S$ and a map $F\colon\cM \ra \cM/S$ such
that the following holds:
\begin{enumerate}[(a)]
\item $F^{cof}$ takes maps in $S$ to weak-equivalences, and
\item $F$ is initial among maps satisfying (a).
\end{enumerate} 
(Recall that $F^{cof}$---which we call the left derived functor of
$F$---denotes any functor obtained by pre-composing $F$ with a
cofibrant-replacement functor.)
\end{defn}

Unfortunately, $S$-localization need not always exist.  The questions
of when they exist and what they might look like can be very hairy, but
there are certain classes of `nice' model categories where the
situation is well under control.  We regard the process of
localization as a way of `imposing relations' in a model category,
hence the notation $\cM/S$---other authors have used $S^{-1}\cM$ or
$L_S\cM$ for the same concept.  

\smallskip

Bousfield \cite{Bo} was the first to give a systematic approach to
what $S$-localizations might look like, and we now recall this.

\begin{defn}\fix
\begin{enumerate}[(a)]
\item An \dfn{$\mdfn{S}$-local} object of $\cM$ is a fibrant object
$X$ such that for every map $A \ra B$ in $S$, the
induced map of homotopy function complexes $\mM(B,X) \ra \mM(A,X)$ is
a weak equivalence of simplicial sets.
\item An \dfn{$\mdfn{S}$-local equivalence} is a map $A\ra B$ such
that $\mM(B,X) \ra \mM(A,X)$ is a weak equivalence for every $S$-local
object $X$.
\end{enumerate}
\end{defn}

\noindent
In words, the $S$-local objects are the ones which see every map in
$S$ as if it were a weak equivalence.  The $S$-local equivalences are
those maps which are seen as weak equivalences by every $S$-local
object.  The idea is that the $S$-local equivalences are the maps
which are {\it forced\/} into being weak equivalences as soon as we
expand our notion of weak equivalence to include the maps in $S$.

\begin{defn}
A \dfn{Bousfield $\mdfn{S}$-localization} of $\cM$ is a model 
category $\cM/S$ with the properties that
\begin{enumerate}[(a)]
\item The underlying category of $\cM/S$ is $\cM$;
\item The cofibrations in $\cM/S$ are the same as those in $\cM$;
\item The weak equivalences in $\cM/S$ are the $S$-local equivalences. 
\end{enumerate}
\end{defn}

Hirschhorn has proven that Bousfield $S$-localizations are indeed
$S$-localizations as defined above.  Bousfield localizations also need
not always exist; if they do exist, they are clearly unique.  The
fibrant objects in $\cM/S$ will be precisely the $S$-local objects,
but the fibrations may be somewhat mysterious.  From now on whenever
we speak of localizations we will always mean Bousfield localizations.

There are two main classes of model categories where localizations are
always known to exist (for any set of maps $S$).  These are the left
proper, cellular model categories of Hirschhorn \cite{H}, and the left
proper, combinatorial model categories of Smith \cite{Sm}. We will not
recall the definitions of these classes here, but suffice it to say
that the model categories $U\cC$ belong to both of them, and so we are
free to localize.  In general, the model categories $U\cC$ are about
as nice as one could possibly want.

\subsection{Basic examples of model categories \mdfn{$U\cC/S$}}
\label{se:locex}

\begin{example} \label{ex:udelta}
In Example~\ref{ex:upt} we saw that the homotopy theory of topological
spaces was the universal homotopy theory on a point, but that this was
almost a tautological statement.  A more interesting example can be
obtained as follows: The way we usually think of simplicial sets is as
objects formally built from the basic simplices, so let us look at
$U\Delta$, the universal homotopy theory built from $\Delta$.  The
obvious map $\Delta \ra \Top$ gives rise to a Quillen pair
$\re:U\Delta \adjoint Top:\sing$, but this is not a Quillen
equivalence.  The first problem one encounters is that there is
nothing in $U\Delta$ saying that the objects $\del{n}$ are
contractible.  In fact this turns out to be the only problem.  If we
localize $U\Delta$ at the set of maps $S=\{\Delta^n \ra *\}$, then our
Quillen functors descend to a pair
\[ \re: U\Delta/\{\Delta^n\ra *\} \adjoint \Top :\sing.\]
It can be seen that this is now a Quillen equivalence---this can be
deduced from \cite[Proposition 5.2]{D1}, but in fact it was the present
observation which inspired that result.  So the
homotopy theory of simplicial sets is the universal homotopy theory
built from $\Delta$ in which the $\Delta^n$'s are contractible.
\end{example}

\begin{example}[Gamma-spaces]
In this example we need the observation that everything we've done
with universal model categories can be duplicated in a {\it pointed\/}
context. Namely, every small category $\cC$ gives rise to a universal
pointed model category built from $\cC$, denoted $U_*\cC$.  Instead of
using presheaves of simplicial sets one uses presheaves of pointed
simplicial sets, and all the same results work with identical
proofs.

Now let $\Spectra$ denote your favorite model category of
spectra---for convenience we'll choose Bousfield-Friedlander spectra---
and let $S$ denote the sphere spectrum.  Let $\cC$ be the subcategory
whose objects are $S$, $S\times S$, $S\times S\times S$, etc., and
whose morphisms are generated by the `obvious' maps one can write down:
e.g., projections $p_1,p_2\colon S\times S \ra S$, inclusions into a factor 
$i_1,i_2\colon S \ra S\times S$, diagonal maps $S \ra S\times S$, etc.   

Now $\cC$ is almost the same as the category called $\Gamma$ in
\cite{BF}---the only difference is that $\Gamma$ contains an extra
object corresponding to the trivial spectrum $*$.  In any case the
inclusion $\cC \inc \Spectra$ extends to a Quillen pair $\re:U_*\cC
\adjoint \Spectra:\sing$, and the category $U_*\cC$ is isomorphic to
the category of $\Gamma$-spaces as defined in \cite{BF}.  The
realization and singular functors are what Segal calls $B$ and $A$,
respectively.  These functors are clearly not a Quillen
equivalence, but let us see if we can somehow turn them into one.

Let $S^{\times k}$ denote the representable object in $U_*\cC$
corresponding to $S\times\cdots\times S$ ($k$ times).  There are
obvious maps $S^{\times k}\Wedge S^{\times l} \ra S^{\times(k+l)}$,
restricting to the inclusions on each wedge-summand---these maps
certainly become equivalences after applying $\re$.  Consider also the
`shearing map' $sh\colon S^{\times 1}\Wedge S^{\times 1} \ra S^{\times
2}$ which on the first wedge-summand is the inclusion $i_1$ and on the
second wedge-summand is the diagonal map.  This map becomes a weak
equivalence under realization as well.  If $W$ denotes the set of all
these maps, then after localizing at $W$ our Quillen pair descends to
give $U_*\cC/W \adjoint \Spectra$.

The model category $U_*\cC/W$ turns out to be precisely one of the
well-known model structures for the category of $\Gamma$-spaces: it is
the one used by Schwede \cite{Sch}, and the identification of the
appropriate maps to localize is implicit in that paper.  The fibrant
objects can be seen to be the `very special' $\Gamma$-spaces (see
\cite[bottom paragraph on Page 349]{Sch} for an argument).  Of course
it's still not quite true that $U_*\cC/W \adjoint \Spectra$ is a
Quillen equivalence: the image of the realization functor consists
only of spectra which can be built from finite products of spheres,
which up to homotopy are the {\it connective\/} spectra.  But it's
well-known that this is the only problem, and that $\Gamma$-spaces
model the homotopy theory of connective spectra.

To summarize: if one starts with a `formal sphere object' $S$
and its finite products $S^{\times k}$, builds the universal pointed
homotopy theory determined by these, and imposes the relations 
\[ S^{\times k}\Wedge S^{\times l} \we S^{\times(k+l)}, \qquad
 sh:S^{\times 1}\Wedge S^{\times 1} \we S^{\times 2}
\]
then one recovers the homotopy theory of connective spectra.  In the
language of section 6 this is a {\it presentation\/} for that
homotopy theory.
\end{example}

\subsection{Further applications}\label{se:cofibmc}

We conclude this section with a result suggesting that model
categories of the form $U\cC/S$ behave something like the {\it
cofibrant objects\/} among model categories.  For another result along
these lines, see Corollary~\ref{co:zigzag}.

\begin{defn} Let $L_1,L_2\colon \cM\ra \cN$ be two maps between 
model categories.
\begin{enumerate}[(a)]  
\item A \dfn{Quillen homotopy} between the maps $L_1$ and $L_2$ is a
natural transformation $L_1 \ra L_2$ with the property that $L_1 X\ra
L_2X$ is a weak equivalence whenever $X$ is cofibrant.
\item As expected, two maps
are \dfn{Quillen-homotopic} if they can be connected by a zig-zag of
Quillen homotopies.
\end{enumerate}
\end{defn}

\begin{prop}
\label{pr:lift}
Let $P\colon\cM\ra \cN$ be a Quillen equivalence of model categories,
and let $F:U\cC/S \ra \cN$ be any map.  Then there is a map $l\colon
U\cC/S \ra \cM$ such that the composite $Pl$ is Quillen-homotopic to
$F$.
\end{prop}

\begin{proof}
See Section~\ref{ad:lifts}.
\end{proof}


\section{Presentations for model categories}

In this section we consider model categories which can be
obtained---up to Quillen equivalence---by starting with a universal
model category $U\cC$ and then localizing at some set of maps $S$. We
refer to these as model categories with {\it presentations\/}, since
the category $\cC$ can be thought of as a category of `generators',
and the set $S$ a collection of `relations'. 

\bigskip

We begin with the basic definition:

\begin{defn}
Let $\cM$ be a model category.  A \mdfn{small presentation} of $\cM$
consists of the following data:
\begin{enumerate}[(1)] 
\item a small category $\cC$, 
\item a choice of Quillen pair $\re:U\cC \adjoint \cM:\sing$,
\item a set of maps $S$ in $U\cC$, 
\end{enumerate}
and we require the properties that
\begin{enumerate}[(a)]
\item The left derived functor of $\re$ takes maps in $S$ to weak
equivalences;
\item The induced Quillen pair $U\cC/S \adjoint \cM$ is a Quillen
equivalence.  
\end{enumerate}
\end{defn}

\noindent
One may think of a small presentation as giving `generators' and
`relations' for the model category $\cM$---see (\ref{se:locex}) in the
preceding section for some examples.  It is not true that every model
category will have a small presentation, but many examples of interest
do.  In fact there is a very large class called the {\it combinatorial
model categories\/} which have been introduced by Jeff Smith, and such
model categories turn out to have small presentations.  Combinatorial
model categories include essentially any model category of algebraic
origin, as well as any model category built-up in some way from
simplicial sets.  We recall the basic definition:

\begin{defn}
A model category $\cM$ is called \dfn{combinatorial} if it is
cofibrantly-generated and the underlying category is locally
presentable.  The latter means that there is a regular cardinal
$\lambda$ and a set of objects $\cA$ in $\cM$ such that
\begin{enumerate}[(i)]
\item Every object in $\cA$ is small with respect to $\lambda$-filtered
colimits, and
\item Every object of $\cM$ can be expressed as a $\lambda$-filtered
colimit of elements of $\cA$. 
\end{enumerate} 
\end{defn}

We can now state the main theorem of this section:

\begin{thm}
\label{th:present}
Any combinatorial model category has a small presentation.
\end{thm}

For background on locally presentable categories we refer the reader
to \cite[Section 1.B]{AR}.  It is a standard result from category
theory that any locally presentable category is equivalent to a full,
reflective subcategory of a category of diagrams $\Set^{\cA}$ (where
$\cA$ is some small category)---see \cite[Prop. 1.46]{AR}.  The
above theorem is the homotopy-theoretic analog of this result.  The
reflecting functor corresponds to the fibrant-replacement functor for
the localized model category.

The proof of Theorem~\ref{th:present} is too involved to give here,
but can be found in the companion paper \cite{D3}.
Here we can at least note two immediate corollaries.  It
was proven in \cite{D1} that any left proper, combinatorial model
category is Quillen equivalent to a simplicial one.  Using the above
theorem we can give a completely different proof of this result, 
and in fact we do slightly better in that we eliminate the left
properness assumption:

\begin{cor}
\label{co:simp}
Any combinatorial model category is Quillen equivalent to one which is
both simplicial and left proper.
\end{cor}

\begin{proof}
The point is that the model categories $U\cC$ are simplicial and left
proper, and these properties are inherited by the localizations $U\cC/S$.
\end{proof}

The second corollary is another instance of the `cofibrant-like'
behavior of the model categories $U\cC/S$---we offer it mainly as an
intriguing curiosity:

\begin{cor} 
\label{co:zigzag}
Suppose one has a zig-zag of Quillen equivalences
\[ \cM_1 \we \cM_2 \bwe \cM_3 \we \cdots \bwe \cM_n \]
in which $\cM_1$ is a combinatorial model category.
Then there is a combinatorial model category $\cN$ and a simple
zig-zag of Quillen equivalences
\[ \cM_1 \bwe \cN \we \cM_n.\]
In fact, $\cN$ may be taken to be of the form $U\cC/S$ where both
$\cC$ and $S$ are small.
\end{cor}

\begin{proof}
One simply chooses a presentation $U\cC/S \we M_1$ and then uses
Proposition~\ref{pr:lift} to lift this map across the Quillen
equivalences.
\end{proof}


\section{Applications to sheaf theory}

Over the years several people have realized that one can construct
model categories which serve as natural settings for `homotopical'
generalizations of sheaf cohomology \cite{BG,Jo,J2}.  What we mean is
that sheaf cohomology appears in these settings as homotopy
classes of maps to certain abelian group objects, but one is allowed
to consider maps to non-additive objects as well.  This `homotopical
sheaf theory' has been very important in applications to algebraic
$K$-theory \cite{Th,J1}, and recently to motivic homotopy theory
\cite{MV}.  In this section we explain a very direct way for
recovering the same homotopy theory via our universal constructions
(\ref{pr:sheaf}).

\bigskip

Recall that a Grothendieck site is a small category $\cC$ equipped
with finite limits, together with a {\it topology}: a collection of
families $\{U_\alpha \ra X\}$ called \dfn{covering families}, which
are required to satisfy various reasonable properties \cite{Ar}.
(There is also a more general approach involving {\it covering
sieves}, which we have foregone only for ease of presentation).  The
prototype for all Grothendieck sites is the category of topological
spaces (contained in a certain universe, say) where the covering
families are just the usual open covers.

If $f\colon E \ra B$ is a map between presheaves, where both $E$ and
$B$ are coproducts of representables, one says that $f$ is a
\dfn{cover} if it has the following property: for any map $rX \ra B$,
there is a covering family $\{U_\alpha \ra X\}$ for which the
compositions $rU_\alpha \ra rX \ra B$ lift through $f$.

\begin{defn}
Let $X\in \cC$ and suppose that $U_*$ is a simplicial presheaf with a
map $U_* \ra rX$.  This map is called
a \dfn{hypercover} of $X$ if
\begin{enumerate}[(i)]
\item Each $U_n$ is a coproduct of representables, 
\item $U_0 \ra rX$ is a cover, and
\item For every integer $n\geq 1$, the component of $U^{\del n} \ra
U^{\bdd n}$ in degree $0$ is a cover.
\end{enumerate}
\end{defn}

This definition is not particularly enlightening, but it's easy to
explain. The easiest examples of hypercovers are the \Cech covers,
which have the form
\[ \xymatrix{  
\cdots \coprod U_{\alpha\beta\gamma} \ar[r]<0.6ex>\ar[r]\ar[r]<-0.6ex> 
    &\coprod U_{\alpha\beta} \ar[r]<0.5ex>\ar[r]<-0.5ex>
    &\coprod U_\alpha \ar[r]
    & X
}
\]
for some chosen covering family $\{U_\alpha \ra X\}$.  Here
$U_{\alpha\beta}=U_\alpha \times_X U_\beta$, etc.  The \Cech covers
are the hypercovers in which the maps $U^{\del n} \ra U^{\bdd n}$ are
{\it isomorphisms} in degree $0$.  In a general hypercover one takes
the iterated fibred-products at each level but then is allowed to
refine that object further, by taking a cover of it.  We refer the
reader to \cite[Section 8]{AM} for further discussion of hypercovers.

For the category of topological spaces one has the following very
useful property: if $U_* \ra X$ is any hypercover of the space
$X$---where in this context we now consider the $U_n$'s as spaces, not
representable presheaves---then the natural map
\[ \hocolim U_* \ra X \]
is a weak equivalence.  (This is not that difficult to prove: if one
has a homotopy element $\bdd{n} \ra X$ then by subdividing the domain
enough times one can gradually lift the map up through the
hypercover).  Based on this observation, it is natural to make the
following construction for any Grothendieck site: 

\begin{defn} Suppose $\cC$ is a Grothendieck site with topology $\cT$. 
Then \mdfn{$U\cC/{\cT}$} denotes the model category obtained by forming
the universal model category $U\cC$ and then localizing at the set of
maps $\{\hocolim U_* \ra X\}$, where $X$ runs through all objects of
$\cC$ and $U_*$ runs through all hypercovers of $X$.
\end{defn}
In words, we have freely added homotopy colimits to $\cC$ and then
imposed relations telling us that any object $X$ may be homotopically
decomposed by taking covers.  Of course sheaf theory is, in the end,
precisely this study of how objects decompose in terms of covers.

\smallskip

In \cite{J2} Jardine introduced a model structure on simplicial
presheaves $\sPre(\cC)$ in which cofibrations are monomorphisms and
weak equivalences are maps inducing isomorphisms on sheaves of
homotopy groups.  We will denote this model category by
$\sPre(\cC)_{Jardine}$.  Since one has the obvious functor $r\colon
\cC \ra \sPre(\cC)$ sending an object to the corresponding
representable presheaf, our general machinery can be seen to give a
map $U\cC/\cT \ra \sPre(\cC)_{Jardine}$ (one must of course check that
the maps we are localizing are weak equivalences in Jardine's sense,
but this is easy).  The essence of the following proposition could
almost be considered folklore---a proof can be found in \cite{D2}:

\begin{prop}
\label{pr:sheaf}
The above map $U\cC/\cT \ra \sPre(\cC)_{Jardine}$ is a Quillen
equivalence.
\end{prop}

\begin{remark}
The model categories $U\cC/\cT$ and $\sPre(\cC)_{Jardine}$ are of
course not that different: they share the same underlying category and
(it turns out) the same weak equivalences, but the notions of
cofibration and fibration differ.  These two different model
structures can already be seen at the level of $U\cC$, before we
localize: in this paper we have consistently used the Bousfield-Kan
model structure, in which the fibrations and weak equivalences are
detected objectwise, but there is also a \dfn{Heller model structure}
\cite{He} in which the {\it cofibrations} and weak equivalences are
detected objectwise.  The Heller structure doesn't seem to enjoy any
kind of universal property, however.
\end{remark}

It is sometimes considered more `natural' to work with simplicial {\it
sheaves} than with simplicial presheaves, although they give rise to
the same homotopy theory---this was what Joyal \cite{Jo} originally
did, and simplicial sheaves were also used in \cite{MV}.  But from the
viewpoint of universal model categories simplicial presheaves are very
natural.  By working with sheaves one allows oneself to recover any
object as a colimit of the objects in a covering, but if you're doing
homotopy theory and only care about {\it homotopy} colimits then
working in the category of sheaves is not so important.  

For more on the rich subject of `homotopical' sheaf theory we refer
the reader to the papers of Jardine \cite{J1,J2} together with
\cite{Be,D2,Th}.

\section{Applications to the homotopy theory of schemes}

Fix a field $k$.  Morel and Voevodsky \cite{MV} have shown that
studying the algebraic $K$-theory and motivic cohomology of smooth
$k$-schemes is part of a much larger subject which they call the {\it
$\A^1$-homotopy theory\/} of such schemes.  They have produced various
Quillen equivalent model categories representing this homotopy theory.
In this section we describe how their procedures relate to our
framework of universal model categories.

\bigskip

Let $\Smk$ denote the category of smooth schemes of finite type over
$k$.  Let $\cT$ be a Grothendieck topology on this category.  Morel
and Voevodsky consider the category of simplicial sheaves
$\sshv(\Smk)$ on this site, with the model structure of \cite{Jo} in
which
\begin{enumerate}[(i)]
\item The cofibrations are the monomorphisms,
\item The weak equivalences are maps which induce weak equivalences on
all stalks (in the case where the site has enough points, which we will
assume for convenience), and
\item The fibrations are the maps with the appropriate lifting
property.
\end{enumerate}
They then define the associated $\A^1$-local structure as the
localization of this model category with respect to the projections
$X\times \A^1 \ra X$, for all $X\in \Smk$.  We'll use the notation
$\MV_k$ for this localized model category.

\smallskip

The point we would like to make is that we can recover the same
homotopy theory from our methods for universal constructions, and in
fact this is not so far from what Morel and Voevodsky actually do.
Based on what we have learned in this paper, it is natural to
construct a homotopy theory for schemes by taking $\Smk$ and expanding
it into the universal model category $U(\Smk)$ by formally adjoining
homotopy colimits.  We will then impose two types of relations:
\begin{enumerate}[(i)]
\item The homotopy-colimit-type relations coming from the Grothendieck
topology, as we saw in the previous section, and
\item The relations $X\times \A^1 \we X$.
\end{enumerate}
Call the resulting model structure $U(\Smk)_{\A^1}$.
The following proposition is essentially routine:

\begin{prop}
\label{pr:MV}
There is a Quillen equivalence $U(\Smk)_{\A^1} \we \MV_k$.  
\end{prop} 

Before giving the proof we need to recall one useful fact. If $L\colon
M \ra N$ is a map between localizable model categories, and $S$ is a
set of maps in $M$, then there is of course an induced map $M/S \ra
N/(L^{cof}S)$.  The fact we need is that if $L$ was a Quillen
equivalence then the induced map on localizations is also a Quillen
equivalence (see \cite{H} for a proof).

\begin{proof}[Proof of (\ref{pr:MV})]
The Yoneda embedding $\Smk \ra \MV_k$ will extend to a map $U(\Smk) \ra
\MV_k$ (and for convenience we choose the extension induced by the
standard cosimplicial resolution, using the fact that $\MV_k$ is a
simplicial model category).  The relations we are imposing in
$U(\Smk)$ clearly hold in $\MV_k$, and so this map descends to
$U(\Smk)_{\A^1} \ra \MV_k$.  It's easy to check that
the left adjoint is the sheafification functor and the
right adjoint is the inclusion of simplicial sheaves into simplicial
presheaves.  

Perhaps the easiest way to see that this is a Quillen equivalence is to
factor the map into two pieces.  In fact, to start with let's forget
about the $\A^1$-homotopy relations; the map we're considering factors
as follows:
\[ U(\Smk)/\cT \adjoint \sPre(\Smk)_{Jardine} \adjoint
\sshv(\Smk)_{Joyal}. \]
Here $U(\Smk)/\cT$ is the model structure constructed in the last
section, and $\sshv(\Smk)_{Joyal}$ is the model structure of \cite{Jo}
mentioned above.  

The first Quillen pair is an equivalence by
Proposition~\ref{pr:sheaf}.  That the second is a Quillen
equivalence is essentially \cite[Prop. 2.8]{J2}.  By the
above observation these also give Quillen equivalences after we
localize at the maps $X\times \A^1 \ra X$.
\end{proof}

For Grothendieck topologies like the Zariski and Nisnevich topologies
one can get by with a much smaller class of relations than the
hypercovers we used above.  In these cases one only has to consider
the \Cech hypercovers coming from certain two-fold covers $\{U_1, U_2
\ra X\}$.  Of course the more manageable the set of relations is, the
better chance one has of understanding the localized model category.
More information about all this can be found in \cite{MV}.

\begin{remark}
Here is one simple instance in which the model category
$U(\Smk)_{\A^1}$ is more handy than $\MV_k$.  Consider the case
where the field is $\C$, so that one has a functor
\[ \Smk \ra \Top, \quad\qquad X \ra X(\C).\]
The functor sends a scheme $X$ to the topological space of its
complex-valued points.  This map immediately induces a Quillen pair
$U(\SmC) \adjoint \Top$, and since the relations that we are localizing    
hold in $\Top$ the Quillen pair descends to
\[ U(\SmC)_{\A^1} \adjoint \Top. \]
One cannot get a similar Quillen pair when using the Morel-Voevodsky
construction, one only gets an adjoint pair on the homotopy
categories---in essence, the model category $\MV_k$ has too many
cofibrations.  Having an actual Quillen pair can be useful, though.
\end{remark}

Concerning our construction of the model category $U(\Smk)_{\A^1}$,
the natural question is how do we know that we have `enough' relations
to give an interesting homotopy theory?  The prototype for this situation
is the case of topological manifolds, in which case these relations
really do generate the usual homotopy theory of topological spaces:

Let $\Man$ denote the category of all topological manifolds which are
contained in $\R^{\infty}$ (the embedding is not part of the data, it
is just a convenient condition to ensure that we have a {\it small\/}
subcategory of manifolds which contains everything we will be
interested in).  This category has a Grothendieck topology consisting
of the usual open covers.  Consider the model category $U(\Man)_{\R}$
obtained by imposing on $U(\Man)$ the same relations we used in
constructing $U(\Smk)_{\A^1}$ (the analog of $\A^1$ is the manifold $\R$).
Note that the obvious map $\Man\ra \Top$ induces a map of model
categories $U(\Man)_{\R} \ra \Top$.

\begin{prop}
The Quillen pair $U(\Man)_{\R} \adjoint \Top$ is a Quillen equivalence.
\end{prop}

\begin{proof}
We only give a sketch.  The reader can also consult
\cite[Prop. 3.3.3]{MV} for a similar statement.

Consider the subcategory $\pt \inc \Man$ whose unique object is the
one-point manifold.  This inclusion induces a Quillen map $U(\pt) \ra
U(\Man)_{\R}$.  The composition
\[ \sSet=U(\pt) \ra U(\Man)_{\R} \ra \Top \]
is the usual realization/singular functor pair, and is therefore a
Quillen equivalence.  So the homotopy theory of topological spaces is
a retract of that of $U(\Man)_{\R}$, and what we have to show is that
$U(\Man)_{\R}$ doesn't contain anything more.  This is where our
relations come in, because they are enough to allow us to unravel any
manifold into a simplicial set.  If $M$ is a manifold we may choose a
cover $U_\alpha$ whose elements are homeomorphic to open balls in
Euclidean space, hence contractible.  For each intersection $U_\alpha
\cap U_\beta$ we may do the same, and so on for all the multiple
intersections---in this way we build a hypercover for $M$ in which all
the open sets are contractible.  Relation (i) allows us to express $M$
(the object in $U(\Man)_\R$) as a homotopy colimit of these
contractible pieces, and relation (ii) allows us to replace each
contractible piece by a point, up to weak equivalence.  So we find
that any representable object in $U(\Man)_\R$ may be expressed as a
homotopy colimit of points, which of course is just the data in a
simplicial set.  In addition we know that every object of $U(\Man)_\R$
is canonically a homotopy colimit of representables, so it follows
that {\it every} object can be decomposed into just a simplicial set.
\end{proof}

\section{The proofs}

This section contains the more technical proofs that were deferred in
the body of the paper.

\subsection{Section 2: Cofibrant replacement in \mdfn{$U\cC$}}
\label{ad:cofrep}
Our first goal in this section is to prove Lemma~\ref{le:cofrep} and
Proposition~\ref{pr:cofrep}.  We must show that given a simplicial
presheaf $F$, the construction $QF$ is a cofibrant-replacement for $F$
in $U\cC$.  We then prove Proposition~\ref{pr:cofrep2}, which is the
same statement for the construction $\cQ F$.

\medskip

Roughly speaking, a simplicial presheaf $F$ will be said to have `free
degeneracies' if there exist presheaves $N_k$ such that $F$ is
isomorphic to the simplicial presheaf
\[ \xymatrix{ 
\cdots
N_2 \amalg (N_1 \amalg N_1 \amalg N_0) \ar[r]<0.6ex>\ar[r]\ar[r]<-0.6ex>
&N_1 \amalg (N_0) \ar[r]<0.5ex>\ar[r]<-0.5ex>
& N_0.}
\]
Here the terms in parentheses in degree $k$ are called the {\it
degenerate part} of $F_k$, and the idea is that these degenerate parts
are as free as possible.  For instance the degenerate part in degree
$2$ consists of a term corresponding to $s_0(N_1)$, a term
corresponding to $s_1(N_1)$, and a term corresponding to
$s_1s_0(N_0)=s_0s_0(N_0)$, and we are requiring that there be no
overlap between these parts.  The following gives a precise
definition:

\begin{defn}  A simplicial presheaf $F$ has \dfn{free degeneracies}
if there exist sub-presheaves $N_k\inc F_k$ such that the canonical
map 
\[ \coprod_{\sigma} N_\sigma \ra F_k
\]
is an isomorphism: here the variable $\sigma$ ranges over all
surjective maps in $\Delta$ of the form $[k]\ra [n]$, $N_\sigma$
denotes a copy of $N_n$, and the map $N_\sigma \ra F_k$ is the one
induced by $\sigma^*\colon F_n \ra F_k$.  (This is called a
\dfn{splitting} of $F$ in \cite[Def. 8.1]{AM}).
\end{defn}

\begin{lemma}
If $F$ has free degeneracies then $F$ is the colimit of the
maps
\[ \sk_0 F \ra \sk_1 F \ra \sk_2 F \ra \cdots \]
where $\sk_0 F=N_0$ and $\sk_n F$ is defined by a pushout-square
\[ \xymatrix{ N_n \tens \bdd{n} \ar[d]\ar[r] 
                       & \sk_{n-1}F \ar@{.>}[d] \\
              N_n \tens \del{n} \ar@{.>}[r] 
                       & \sk_n F.}
\]
\end{lemma}

\begin{proof}
Left to the reader.
\end{proof}

\begin{cor}
\label{co:cofsplit}
If $F$ has a free degeneracy decomposition in which the $N_k$ are
cofibrant in $U\cC$, then $F$ is itself cofibrant.   
\end{cor}

\begin{proof}
The fact that $N_k$ is cofibrant implies that $N_k\tens \bdd{k} \ra
N_k\tens \del{k}$ is a cofibration, and so the map $\sk_{k-1}F \ra
\sk_k F$ is also a cofibration.  Then $F$ is a sequential colimit of
cofibrations beginning with $\emptyset \cof \sk_0 F$, hence
cofibrant.
\end{proof}

We now prove that if $F$ is a {\it discrete\/} simplicial presheaf
then $\tilde{Q}F$ is a cofibrant replacement for $F$:

\begin{proof}[Proof of Lemma~\ref{le:cofrep}]
First observe that $\tilde{Q}F$ has a free degeneracy decomposition:
we take $N_k$ to be the coproduct
\[ \coprod_{rX_k \ra \cdots rX_0\ra F} (rX_k) \]
in which no map $X_{i+1} \ra X_i$ is an identity map.
Each $N_k$ is a coproduct of representables, hence cofibrant.
So $\tilde{Q}F$ is itself cofibrant by the above corollary.

We must next show that $\tilde{Q}F \ra F$ is a weak equivalence in
$U\cC$---that is, we must show that $(\tilde{Q}F)(X) \ra F(X)$ is a
weak equivalence of simplicial sets, for every $X\in \cC$.  
Let $\cA$ denote the subcategory of $\cC$ consisting of the same
objects but only identity maps. 
Consider the adjoint pair
\[ T:\Set^{\cA^{op}} \adjoint \Set^{\cC^{op}}:U, \]
where $U$ is the restriction functor and $T$ is its left adjoint.
Then $TU$ is a cotriple, and the cotriple resolution
\[ \xymatrix{
\cdots
(TU)^3 F \ar[r]<0.6ex>\ar[r]\ar[r]<-0.6ex>
&(TU)^2 F \ar[r]<0.5ex>\ar[r]<-0.5ex>
&(TU)F \ar[r] 
&F}
\]
can be seen to exactly coincide with $\tilde{Q}F$.  Now of course if
we apply $U$ again then we pick up an extra degeneracy, and the map
$U[(TU)^*F] \ra UF$ is a weak equivalence in $\sSet^{\cA^{op}}$.  
But applying $U$ to a simplicial presheaf gives precisely the
collection of all its values, and so we have that $({\tilde
Q}{F})(X) \ra F(X)$ is a weak equivalence for every $X$.
\end{proof}

Now we move on to handle arbitrary simplicial presheaves:

\begin{proof}[Proof of Proposition~\ref{pr:cofrep}]
One again shows that $QF$ has a free degeneracy decomposition in which
the $N_k$ are coproducts of representables.  This takes a little more
work than for $\tilde{Q}F$, but we will leave it to the reader.  The
fact that $QF$ is cofibrant follows from Corollary~\ref{co:cofsplit}.

To see that $QF \ra F$ is a weak equivalence we consider the
bisimplicial object $Q_{**}F$ whose $n$th row is $\tilde Q(F_n)$, as
well as the `constant' bisimplicial object $F_{**}$ whose $n$th row is
the discrete simplicial presheaf consisting of $F_n$ in every level.
The map $QF\ra F$ is the diagonal of a map $Q_{**}F \ra F_{**}$.  But
we have already shown that $\tilde Q{F_n} \ra F_n$ is a weak
equivalence for every $n$, which says that $Q_{**}F \ra F_{**}$ is a
weak equivalence on each row.  It follows that the map yields a weak
equivalence on the diagonal as well.
\end{proof}

The last thing we must do is prove Proposition~\ref{pr:cofrep2}, which
concerned a different functor $\cQ F$---we are to show that this is
another cofibrant-replacement functor for $U\cC$.  The proof is an
unpleasant calculation of a homotopy colimit.

\begin{proof}[Proof of Proposition~\ref{pr:cofrep2}]
We explained in section (\ref{se:cofrep}) why $\cQ F$ was cofibrant,
therefore the only thing to prove is that the natural map $\cQ F \ra
F$ is a weak equivalence in $U\cC$.  So we need to show that for every
$X \in \cC$ the map $\cQ F(X) \ra F(X)$ is a weak equivalence of
simplicial sets.

For brevity let $I$ denote the category $(\cC\times\Delta\ovcat F)$.
The object $\cQ F$ is the homotopy colimit of the diagram $I\ra U\cC$
which sends $[A\times [n], rA\tens\del{n} \ra F]$ to $rA\tens\del{n}$.
Because the simplicial structure in $U\cC$ is the objectwise
structure, homotopy colimits are also computed objectwise.  This says
that $\cQ F(X)$ is equal to the homotopy colimit of the diagram
$\cD\colon I\ra \sSet$ sending $[A\times [n], rA\tens\del{n}\ra F]$ to
$(rA\tens\del{n})(X)$.  This latter object may be identified with
$rA(X) \tens \del{n}$, which is $\cC(X,A)\tens
\del{n}$---it is a coproduct of copies of $\del{n}$, one for each map
$X\ra A$.

Consider the functor $\Theta\colon I \ra \Set$ which sends $[A\times [n],
rA\tens \del{n}\ra F]$ to the set $\cC(X,A)$.  From this functor we
may form its Grothendieck construction $\Gr \Theta$: this is the
category whose objects are pairs $(i,\sigma)$ where $i\in I$ and
$\sigma \in \Theta(i)$, and a map $(i,\sigma) \ra (j,\alpha)$ is a map
$i\ra j$ in $I$ such that $(\Theta i)(\sigma)=\alpha$.  An object of
$\Gr \Theta$ corresponds to the data $[A\times [n], rA\tens\del{n} \ra F,
X\ra A]$, so define a functor $\cE\colon \Gr \Theta \ra \sSet$ 
which sends this object to the simplicial set $\del{n}$.

Thomason has a theorem about homotopy colimits over Grothendieck
constructions \cite[Cor. 24.6]{CS}, and in our situation it gives us a
weak equivalence
\[ \hocolim_{\Gr \Theta} \cE \we 
         \hocolim_{i \in I} 
              \bigl [\hocolim_{\sigma\in \Theta(i)} \cE(i,\sigma) 
              \bigr ].
\]
If $i\in I$ corresponds to the data $[A\times [n], rA\tens\del{n} \ra
F]$, then the homotopy colimit inside the brackets is just a coproduct
of copies of $\del{n}$, one for each element of $\Theta(i)=\cC(X,A)$.
In this way the double homotopy colimit on the right is readily
identified with $\hocolim_I \cD$, and we have already seen that this
is $\cQ F(X)$.

Now consider the category $\Delta(X,F)$, defined so that the objects
consist of the data $[[n], rX\tens \del{n} \ra F]$---this is equal to
the {\it category of simplices\/} of the simplicial set $F(X)$
(defined in \cite[p. 75]{Ho}, for instance).  We again let $\cE \colon
\Delta(X,F) \ra \sSet$ denote the diagram which sends $[[n],
rX\tens\del{n} \ra F]$ to $\del{n}$.  The colimit of this diagram is
just $F(X)$, and the natural map $\hocolim_{\Delta(X,F)} \cE \ra F(X)$
is a weak equivalence of simplicial sets.

There is a functor $\Delta(X,F) \ra \Gr \Theta$ which sends
\[ [[n], rX\tens \del{n} \ra F] \mapsto 
[X\times [n], X\tens \del{n} \ra F, \Id\colon X\ra X], \]
and this induces a map of homotopy colimits
$\hocolim_{\Delta(X,F)} \cE \ra \hocolim_{\Gr \Theta} \cE$.  The map
of categories has a retraction which is easily checked to be
homotopy-cofinal, so it follows that the map of homotopy colimits is a
weak equivalence.

All-in-all what we have is the following diagram:
\[ \xymatrix{ \hocolim_{\Delta(X,F)} \cE \ar[r]\ar[rrd] & \hocolim_{\Gr
\Theta} \cE \ar[r] & \hocolim_I \cD \ar[d] \\
&&F(X).
}
\]
We have shown that every map is a weak equivalence except the vertical
one, but then the vertical map must be one as well.  This is the statement
that $\cQ F(X) \ra F(X)$ is a weak equivalence, which was our goal.
\end{proof}

\subsection{Section 3: Cosimplicial resolutions and maps from universal 
model categories}
\label{ad:fact}
In this section we prove Proposition~\ref{pr:fact}, which said that
extending a map $\gamma\colon\cC \ra \cM$ to the universal model
category $U\cC$ was equivalent to giving a cosimplicial resolution on
$\gamma$.

\medskip

To begin with we will need some machinery.
If $K\in \sSet$ and $X^*\in c\cM$ one can define a tensor product
$X\tens K\in \cM$ (see \cite[Prop. 3.1.5]{Ho}).  Start with some general
notation: For a set $S$ and an object $W \in \cM$, let $W\cdot S$
denote a coproduct of copies of $W$, one for each element of $S$.
Then $X\tens K$ can be defined as a coend:
\[ X\tens K=\coeq \Biggl[ \coprod_{[k]\ra[m]} X_k \cdot K_m \dbra 
                      \coprod_n X_n \cdot K_n \Biggr ].\]
This construction has the adjointness property that
\[ \cM(X\tens K,W)\cong \sSet(K,\cM(X^*,W)) \]
where $\cM(X^*,W)$ is the simplicial set whose $n$-simplices are the
hom-set $\cM(X^n,W)$.   

Now if we have diagrams $\Gamma\colon\cC \ra c\cM$ and
$F\colon\cC^{op} \ra \sSet$ then we can again form a coend
\[ \Gamma \tens_{\cC} F
          =\coeq \Biggl[ \coprod_{a\ra b} \Gamma(a) \tens F(b) \dbra 
                      \coprod_{c\in \cC} \Gamma(c) \tens F(c) \Biggr
                      ].
\]
For this construction we have that
\begin{equation}
\label{eq:adjnt}
 \cM(\Gamma\tens_{\cC} F,W)=\sSet^{\cC^{op}}(F,\cM(\Gamma,W)) 
\end{equation}
where $\cM(\Gamma,W)$ is the simplicial presheaf defined by
$ c \mapsto \cM(\Gamma^*c,W)$.

The above is all that's necessary to prove our result:

\begin{proof}[Proof of Proposition~\ref{pr:fact}]
Suppose we have a factorization of $\gamma\colon\cC \ra \cM$ through
$U\cC$: so we have a Quillen pair $\re:U\cC \adjoint \cM:\sing$ and a
natural weak equivalence $\re(rX)\we \gamma(X)$.  Then for each $X\in
\cC$ we get a cosimplicial resolution of $\gamma X$ by taking
$\Gamma(X)$ to be
\[ [n] \mapsto \re(rX\tens \del{n}).\]
This is clearly functorial in $X$, and so gives a resolution of $\gamma$.

Conversely, suppose we start with a resolution $\Gamma\colon\cC\ra c\cM$
for the functor $\gamma$.  Define the functors $\re:U\cC \ra \cM$ and
$\sing: \cM \ra U\cC$ by the formulas
\[ \re(F)=\Gamma\tens_\cC F, \qquad   
             \sing(X)=[c\mapsto \cM(\Gamma^*(c),X)].
\] 
(\ref{eq:adjnt}) says that these are an adjoint pair.

To see that these are a Quillen pair we will check that $\sing$
preserves fibrations and trivial fibrations.  For this we need to know
that if $A^*$ is a cosimplicial resolution and $X\ra Y$ is a fibration
(resp. trivial fibration) then $\cM(A^*,X) \ra \cM(A^*,Y)$ is a
fibration (resp. trivial fibration) of simplicial sets.  But this is
\cite[Cor. 5.1.4]{Ho}.

The last thing is to give a natural weak equivalence $\re(rX) \ra
\gamma(X)$.  But $\re(rX)$ is isomorphic to the object
of $\Gamma(X)$ in degree $0$, and our cosimplicial resolution came
with a weak equivalence from this object to $\gamma(X)$.  So we're
done.

Checking the equivalence of categories $\Fact(\gamma) \he
\coRes(\gamma)$ is fairly routine at this point: we have given the
functors in either direction.
\end{proof}


\subsection{Section 5: Lifting maps from the model categories
\mdfn{$U\cC/S$}}
\label{ad:lifts}
Here we fill in the proof of Proposition~\ref{pr:lift}.  We must show
that a map from a model category $U\cC/S$ may be lifted, up to
homotopy, across a Quillen equivalence.

\medskip

It will be useful to isolate the following lemma:

\begin{lemma}
\label{le:quilho}
Let $\cM$ be a model category, and let $\gamma_1,\gamma_2\colon\cC \ra
\cM$ be two functors whose images lie in the cofibrant objects.
Suppose there is a natural weak equivalence $\gamma_1 \we \gamma_2$.
Then any two extensions $L_1,L_2\colon U\cC \ra \cM$ of $\gamma_1$ and
$\gamma_2$ are Quillen-homotopic.
\end{lemma}

\begin{proof}
Recall that there exists an equivalence between maps of model
categories $U\cC \ra \cM$ and the following data:
\begin{itemize}
\item A functor $\gamma\colon\cC \ra \cM$ whose image lies in the cofibrant
objects, and
\item A cosimplicial resolution on $\gamma$.
\end{itemize}
Giving a Quillen homotopy between two maps $L_1,L_2\colon U\cC \ra \cM$
corresponds to giving a natural weak equivalence $\gamma_1 \ra
\gamma_2$ and a lifting of this to a natural weak equivalence between
the resolutions.  Using these facts, proving the lemma is just a matter
of getting zig-zags between the resolutions.  But this is standard---see
\cite{H}.
\end{proof}

\begin{proof}[Proof of Proposition~\ref{pr:lift}]
Let $\tilde F$ be the composite $U\cC \ra U\cC/S \ra \cN$.  We will
begin by lifting $\tilde F$, and this can be accomplished just by
lifting $\gamma\colon \cC \ra U\cC \ra \cN$.

Define $\epsilon:\cC \ra \cM$ by
\[ \epsilon(X)=[Q^{fib}(\gamma(X))]^{cof}\]
where $Q$ is the right-adjoint to $P$.  We may extend $\epsilon$ to a
map $\tilde l:U\cC \ra \cM$.

\smallskip

\noindent{\it Claim 1:\/} The composite $P\tilde{l}$ is
Quillen-homotopic to $\tilde F$.

\smallskip\noindent
To see this, observe that there are natural weak equivalences
\[ P\epsilon(X) \we [\gamma X]^{fib} \bwe \gamma(X).\]
Since $P\tilde l$ is an extension of $P\epsilon$ and $F$ is an
extension of $\gamma$, the claim follows directly from the above
lemma.

\medskip

\noindent{\it Claim 2:\/} The map $\tilde{l}$ takes elements of $S$ to weak
equivalences.

\smallskip\noindent For this, note that by hypothesis the derived
functor of $\tilde F$ takes elements of $S$ to weak equivalences.  The
same must be true for $P\tilde l$, since $P\tilde l$ is homotopic to
$\tilde F$ (maps which are Quillen-homotopic will have isomorphic
derived functors on the homotopy categories).  But $P$ was a Quillen
equivalence, and so the derived functor of $\tilde l$ must also take
elements of $S$ to weak equivalences.

\medskip

From Claim 2 it follows that $\tilde l$ descends to a map $l:U\cC/S
\ra \cM$.  The fact that $Pl$ is homotopic to $F$ is just a restatement
of Claim 1.
\end{proof}


\bibliographystyle{amsalpha}

\end{document}